\documentclass[11pt]{amsart}
\usepackage{amsmath,amsthm,amsfonts,amssymb}
\usepackage{pictex}
\usepackage{youngtab}
\usepackage{young}
 
\newtheorem{thm}[equation]{Theorem} 
\newtheorem{cor}[equation]{Corollary}
\newtheorem{lem}[equation]{Lemma}

\theoremstyle{definition}  
\newtheorem{defn}[equation]{Definition} 

\theoremstyle{remark}
\newtheorem{rem}[equation]{Remark} 
\numberwithin{equation}{section}

\newcommand{\z}{\mathbb{Z}}

\begin{document}
\title[Constructing all irreducible Specht modules in a block]{Constructing all irreducible Specht modules in a block of the symmetric group}

\author{James P. Cossey}
  \address{Department of Mathematics\\
           University of Arizona\\
           Tucson, AZ  85721  USA}
  \email{cossey@math.arizona.edu}
\author{Matthew Ondrus}
  \address{Department of Mathematics\\
           University of Arizona\\
           Tucson, AZ  85721  USA}
  \email{ondrus@math.arizona.edu}
\author{C. Ryan Vinroot} 
  \address{Department of Mathematics\\
           University of Arizona\\
           Tucson, AZ  85721  USA}
  \email{vinroot@math.arizona.edu}

\begin{abstract}
For any prime $p$, we construct, and simultaneously count, all of the complex Specht modules in a given $p$-block of the symmetric group which remain irreducible when reduced modulo $p$.  We call the Specht modules with this property $p$-irreducible modules.  Recently Fayers has proven a conjecture of James and Mathas that provides a characterization of the partitions that correspond to the $p$-irreducible modules.  In this paper we present a method for decomposing the partitions corresponding to $p$-irreducible modules, and we use this decomposition to construct and count all of the partitions corresponding to $p$-irreducible Specht modules in a given block.
\end{abstract}

\maketitle

\section{Introduction}\label{I}

Recently, Fayers \cite{F1} has proven a conjecture which gives a condition on a partition, $\lambda$, so that the Specht module  $S^\lambda_{\mathbb F_p}$ defined over $\mathbb F_p$, for $p > 2$, is irreducible.  Equivalently, this is a condition on $\lambda$ so that the Specht module $S^{\lambda}$ defined over the complex field remains irreducible when reducing it modulo $p$.  In this paper, we begin by giving a decomposition of a partition $\lambda$ which satisfies this condition.  We then use this decomposition to describe which of these partitions correspond to representations in a given $p$-block of the symmetric group.

In order to state our main results, we must first clarify several definitions.  We use standard notations and definitions for partitions and representations of the symmetric group, which are stated in Section 2.  For any prime $p$, we call a Specht module, $S^{\lambda}$, \emph{$p$-irreducible} if $S^{\lambda}$ is defined over the complex field and remains irreducible when reduced modulo $p$, or equivalently, if $S^{\lambda}_{\mathbb F_p}$ is irreducible.  In Section 2, we state Fayers' theorem (Theorem \ref{fayers}), which gives a condition on $\lambda$, for $p>2$, which describes exactly when $S^{\lambda}$ is $p$-irreducible.  We will say that the partitions which meet this condition are \emph{$p$-irreducible} (even if $p=2$).  The case $p=2$ is slightly different and was described by James and Mathas in \cite{JM}.  We discuss this case in Section \ref{twocase}.  

We recall that by Nakayama's conjecture, two Specht modules of the symmetric group $S_n$, $S^{\lambda}$ and $S^{\mu}$, defined over the complex field, are in the same $p$-block precisely when $\lambda$ and $\mu$ have the same $p$-core, which is the result of removing rims of all $p$-hooks of a partition.  We call a partition \emph{$p$-hook free} when none of its hook-lengths are divisible by $p$.  The $p$-blocks of symmetric groups are therefore represented by $p$-hook free partitions.  If the $p$-block $B$ of $S_n$ is represented by the $p$-hook free partition $\nu$, then the \emph{weight} of $B$ is the number $w$ which satisfies $|\nu| + pw = n$.  For a partition $\nu$, we let $\nu'$ be the conjugate partition and $\ell(\nu)$ to be the number of parts of $\nu$.  Recall that a partition is $p$-regular if none of its parts are repeated $p$ or more times.  We need to define the following feature of a $p$-block.

\begin{defn} \label{residual}
Let $\nu = (\nu_1, \nu_2, \ldots, \nu_l)$ be any partition, where we set $\nu_j = 0$ for $j > l$.  If $\nu_1 - \nu_2 \neq p-1$, let $t = 1$, and if $\nu'_1 - \nu'_2 \neq p-1$ let $b = 1$, and if $\nu$ is the empty partition let $t=b=1$.  Otherwise, we let $t$ and $b$ be the positive integers such that $\nu_i - \nu_{i+1} = p-1$ for $i<t$, but $\nu_t - \nu_{t+1} \neq p-1$, and $\nu'_i - \nu'_{i+1} = p-1$ for $i<b$, but $\nu'_b - \nu'_{b+1} \neq p-1$.  Then we call the ordered pair of positive intergers $(t,b)$ the \emph{$p$-residual} of $\nu$.  If $B$ is a $p$-block of $S_n$ represented by the $p$-hook free partition $\nu$ which has $p$-residual $(t,b)$, then we say the block $B$ has $p$-residual $(t,b)$ as well.
\end{defn}

\noindent
We may now state the main results of the paper, the proofs of which are given in Section \ref{mainthm}.  

\vspace{.15cm}

\noindent \textbf{I.} (Theorem \ref{COVproof}) Let $p > 2$, and let $B$ be a $p$-block of $S_n$, corresponding to the $p$-hook free partition $\nu$, of weight $w$ and $p$-residual $(t,b)$.  Then the number of $p$-irreducible Specht modules in $B$ is equal to the number of ordered pairs of $p$-regular $p$-irreducible partitions $(\alpha, \gamma)$, such that $|\alpha| + |\gamma| = w$, $\ell(\alpha) \leq t$, $\ell(\gamma) \leq b$, and $\ell(\alpha) + \ell(\gamma) \leq t+b-1$ if $t+b = \frac{\ell(\nu)+ \ell(\nu')}{p} + 2$.

\vspace{.15cm}

\noindent \textbf{II.} (Theorem \ref{twocount}) Let $n \neq 4$, and let $B$ be a $2$-block of $S_n$, corresponding to the $2$-hook free partition $\nu$, of weight $w$.  Then the number of $2$-irreducible Specht modules in $B$ is equal to twice the number of $2$-regular $2$-irreducible partitions $\alpha$ of $w$ such that $\ell(\alpha) \leq \ell(\nu)+1$.

\vspace{.15cm}

We note that the proofs of these theorems are constructive, meaning that given any ordered pair of partitions $(\alpha, \gamma)$, or single partition in the case $p=2$, with the properties above, we can construct the corresponding $p$-irreducible module $S^{\lambda}$ in $B$, and vice versa.  We demonstrate this in Section \ref{example} by finding all $5$-irreducible Specht modules in a particular block of $S_{110}$.

It should be noted here that much of this paper could be done in terms of abaci.  For example, our Theorem \ref{construction} and Lemma \ref{shrinking} imply Fayers' characterization of abaci corresponding to $p$-irreducible partitions obtained in \cite[Proposition 2.1]{F1}.  Our arguments deal with partitions directly, rather than the abaci, and provide a direct decomposition and characterization of the $p$-irreducible partitions which can be easily pictured.

\section{Notation and Background}

A {\em partition}, $\lambda$,  of a non-negative integer $n$ is a sequence of {\em parts}, $\lambda = (\lambda_1, \lambda_2, \ldots, \lambda_k)$, where $\lambda_i \geq \lambda_{i+1}$, each $\lambda_i >0$ for $1 \leq i \leq k$, and $\sum_{i=1}^k \lambda_i = n$.  We may set $\lambda_j = 0$ for $j > k$.  The {\em empty partition} is the partition with no nonzero parts.  The {\em length} of $\lambda$ is the number of nonzero parts of $\lambda$, written $\ell(\lambda) = k$.  We may represent a partition by its {\em Young diagram}, which is given by rows of boxes, where the number of boxes in each row corresponds to the parts of $\lambda$, decreasing from top to bottom.  We also use the standard exponential notation to indicate a row that is repeated more than once.  For example, the Young diagram of $\lambda = (7, 3, 2^2, 1)$ is:

$$\yng(7,3,2,2,1)$$

We label the top left corner box of a Young diagram with the coordinates $(1,1)$, and the coordinates of boxes, when moving to the right, increase by $1$ in the second coordinate for each box, and when moving down, increase by $1$ in the first coordinate for each box.  So, for example, the coordinates of the boxes in the above Young diagram, when moving from left to right in the second row from the top, are $(2,1)$, $(2,2)$, and $(2,3)$.  We call each $(i,j)$ in a Young diagram of $\lambda$ a {\em node} of $\lambda$.

The {\em conjugate} of a partition $\lambda$, written $\lambda'$, is the partition obtained when changing the columns of the Young diagram of $\lambda$ into rows.  The parts of $\lambda'$ are given by $\lambda'_i = \big|\{j | \lambda_j \geq i\} \big|$.  The conjugate partition of $(7, 3, 2^2, 1)$, for example, is $(5, 4, 2, 1^3)$, pictured below.

$$\yng(5,4,2,1,1,1)$$

The {\em hook length} at the position $(i,j)$ of $\lambda$, written $h_{\lambda}(i,j)$, is the number of boxes in the Young diagram of $\lambda$ to the right and below $(i,j)$, and the box at $(i,j)$ itself.  That is, $h_{\lambda}(i,j) = \lambda_i - i + \lambda'_j - j + 1$.  Below is the Young diagram of $(7, 3, 2^2, 1)$, where the number in each box is the hook length of that node.
$$\begin{Young}
     11  &  9 &  6 & 4 & 3 & 2 & 1\cr
     6  &  4 & 1     \cr
     4  &  2     \cr
     3 & 1           \cr
     1              \cr
\end{Young}$$

For a positive integer $m$ and a prime $p$, let $v_p(m)$ be the highest power of $p$ which divides $m$.  We say there is a {\em $p$-hook} at $(i,j)$ of a partition $\lambda$ if $v_p(h_{\lambda}(i,j)) > 0$, and that a partition $\lambda$ is {\em $p$-hook free} if at every position $(i,j)$ of $\lambda$ we have $v_p(h_{\lambda}(i,j)) = 0$.  

A partition $\lambda$ is called {\em $p$-regular} if $\lambda$ has no nonzero part which occurs $p$ or more times, and $\lambda$ is called {\em $p$-restricted} if $\lambda'$ is $p$-regular.

Given a partition $\lambda$ of the positive integer $n$, and any field $F$, one may construct an $FS_n$-module called the {\em Specht module}, denoted $S^{\lambda}_F$.  For the construction and basic theory of Specht modules, see, for example \cite[Chapter 7]{JK}.  When $F$ has characteristic $0$, we denote the Specht module by $S^{\lambda}$, which is irreducible, and as $\lambda$ ranges over all partitions of $n$, the modules $S^{\lambda}$ range over all non-isomorphic irreducible $FS_n$-modules.  

If $F$ has characteristic $p>0$, then $S^{\lambda}_F$ may be reducible, but when $\lambda$ is $p$-regular, then $S^{\lambda}_F$ has a uniquely defined quotient module, $D^{\lambda}$, which is irreducible.  As $\lambda$ ranges over all $p$-regular partitions of $n$, the $D^{\lambda}$ range over all distinct irreducible $FS_n$-modules.  We note that when $F$ has characteristic $p>0$, we may always choose a basis so that $S^{\lambda}_F$ may be viewed as being over $\mathbb F_p$  

If $S^{\lambda}$ is defined over a field of characteristic $0$, then since one can choose a basis for $S^{\lambda}$ such that the entries of any matrix corresponding to any group element under this representation are integers, we may reduce $S^{\lambda}$ modulo a prime $p$, and what we get is exactly the Specht module $S^{\lambda}_{\mathbb F_p}$ defined over that characteristic $p$ field.  We say that $S^{\lambda}$ is {\em $p$-irreducible} if it is defined over a characteristic $0$ field, and remains irreducible when reducing it modulo the prime $p$, and this is equivalent to the Specht module $S^{\lambda}_{\mathbb F_p}$ being irreducible.

James and Mathas \cite{M} gave a condition (given in Theorem \ref{fayers} below) on the partition $\lambda$ which they conjectured was necessary and sufficient for $S^{\lambda}$ to be $p$-irreducible for $p > 2$.  This condition was known to be necessary and sufficient for $S^{\lambda}$ to be $p$-irreducible, $p>2$, when $\lambda$ is a $p$-regular or $p$-restricted partition, as proven in \cite{J} and \cite{JMu}.  A significant step in the necessity of this condition in the general case was proven by Lyle \cite{L} in 2003, and completed by Fayers \cite{F2} in 2004.  In 2005, Fayers finished proving this conjecture by proving the sufficiency in \cite{F1}.  The result is as follows.

\begin{thm}[Fayers, 2005] \label{fayers}
Let $p$ be a prime, $p>2$, $n$ a positive integer, and let $\lambda$ be a partition of $n$.  Then the Specht module $S^{\lambda}$ of $S_n$ is $p$-irreducible if and only if $\lambda$ does not have positions $(i,j)$, $(i,y)$, and $(x,j)$ such that 
$$v_p(h_{\lambda}(i,j)) > 0, \; v_p(h_{\lambda}(i,y)) \neq v_p(h_{\lambda}(i,j)), \text{ and } v_p(h_{\lambda}(x,j)) \neq v_p(h_{\lambda}(i,j)).$$
\end{thm}

We thus say that partitions $\lambda$ which meet the condition in Theorem \ref{fayers} are \emph{$p$-irreducible}.  Even for $p=2$, we say that $\lambda$ is $p$-irreducible when the condition in Theorem \ref{fayers} holds, even though this is not equivalent to $S^{\lambda}$ being $2$-irreducible in general, as we see in Section \ref{twocase}.  

\section{The decomposition of $p$-irreducible partitions} \label{decomposition}

Given three partitions $\tau$, $\mu$, and $\beta$, we construct a new partition containing all three through a specific gluing, as we now describe.  

\begin{defn} \label{decomp} Let $\tau$, $\mu$, and $\beta$ be three partitions, where $\tau$ and $\beta$ may be empty, but $\mu$ may only be empty if at least one of $\tau$ or $\beta$ is empty.  If $\tau = (\tau_1, \ldots, \tau_k)$, $\mu = (\mu_1, \ldots, \mu_l)$, and $\beta = (\beta_1, \ldots, \beta_m)$, then define the partition $\oplus(\tau, \mu, \beta)$ to be the partition with parts $(\tau_1 + \mu_1 + \beta_1 -1, \, \tau_2 + \mu_1 + \beta_1 - 1, \ldots, \tau_k + \mu_1 + \beta_1 - 1, \, \mu_2 + \beta_1 - 1, \, \mu_3 + \beta_1 -1, \ldots, \mu_l + \beta_1 - 1, \, \beta_1, \ldots, \beta_m)$.
\end{defn}

One may visualize this gluing in the following way.  We glue the first node in the last row of $\tau$ to the right of the last node of the first row of $\mu$, and glue the last node of the first column of $\mu$ to the top of the first node of the last column of $\beta$.  We then fill in the partition in the natural way to the upper-left corner.  Note that if $\mu$ is empty, then either $\tau$ or $\beta$ is also, and there is no gluing.  This gluing looks like the following picture:

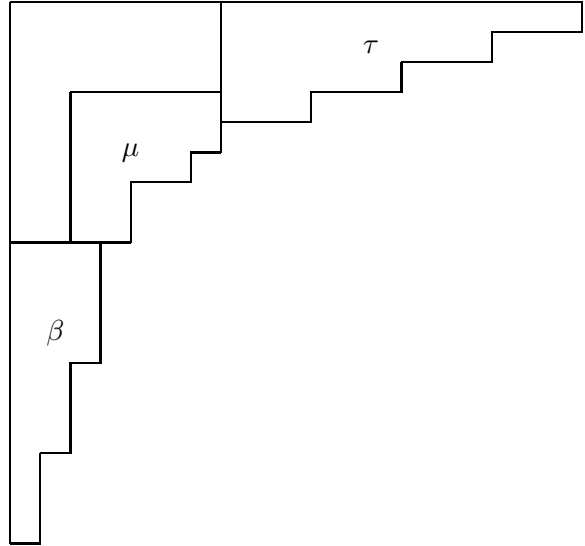
\begin{center}
$\beginpicture
\setcoordinatesystem units <0.4cm,0.4cm>         
\setplotarea x from 0 to 19, y from 0 to 19
\putrule from 0 0 to 1 0
\putrule from 1 3 to 2 3
\putrule from 2 6 to 3 6
\putrule from 0 10 to 3 10

\putrule from 0 0 to 0 10
\putrule from 1 0 to 1 3
\putrule from 2 3 to 2 6
\putrule from 3 6 to 3 10

\putrule from 2 10 to 4 10
\putrule from 4 12 to 6 12
\putrule from 6 13 to 7 13
\putrule from 2 15 to 7 15

\putrule from 2 10 to 2 15
\putrule from 4 10 to 4 12
\putrule from 6 12 to 6 13
\putrule from 7 13 to 7 15

\putrule from 7 14 to 10 14
\putrule from 10 15 to 13 15
\putrule from 13 16 to 16 16
\putrule from 16 17 to 19 17
\putrule from 7 18 to 19 18

\putrule from 7 14 to 7 18
\putrule from 10 14 to 10 15
\putrule from 13 15 to 13 16
\putrule from 16 16 to 16 17
\putrule from 19 17 to 19 18

\putrule from 0 10 to 0 18
\putrule from 0 18 to 7 18

\put{$\tau$} at 12 16.5
\put{$\mu$} at 4 13
\put{$\beta$} at 1.5 7

\endpicture$\\
\end{center}

We are specifically interested in considering $\oplus(\tau, \mu, \beta)$ when $\mu$ is $p$-hook free, and when $\tau$ and $\beta$ are certain types of partitions, which we now define.

\begin{defn}
A \emph{p-top} is a partition $\tau = (\tau_1, \ldots, \tau_m)$ such that $v_p(h_{\tau}(i,1)) > 0$ for $1 \leq i \leq m$.  A \emph{p-irreducible top} is a $p$-top which is also $p$-irreducible.

Likewise, a \emph{p-bottom} is a partition $\beta = (\beta_1, \ldots, \beta_m)$ such that $v_p(h_{\beta}(1,i)) > 0$ for $1 \leq i \leq \beta_1$, and a \emph{p-irreducible bottom} is $p$-bottom which is also $p$-irreducible. 
\end{defn}

\noindent
Note that $\tau$ is a $p$-top or $p$-irreducible top if and only if $\tau'$ is a $p$-bottom or $p$-irreducible bottom, repsectively.  We now give the decomposition of a $p$-irreducible partition.
 
\begin{thm} \label{construction}
If $\lambda$ is a $p$-irreducible partition, then we may decompose $\lambda$ as $\lambda = \oplus(\tau, \mu, \beta)$ for a unique $p$-irreducible top $\tau$, $p$-hook free $\mu$, and $p$-irreducible bottom $\beta$.  Conversely, for any $p$-irreducible top $\tau$, $p$-hook free $\mu$, and $p$-irreducible bottom $\beta$, the partition $\oplus(\tau, \mu, \beta)$ is $p$-irreducible.
\end{thm}

\begin{proof}
We first prove that if $\lambda = \oplus(\tau, \mu, \beta)$ is constructed
as above, then $\lambda$ is $p$-irreducible.  We do this by showing
that if $(i, j)$ is a node of $\lambda$ not contained in $\tau$ or
$\beta$, that is, if $i < \tau'_1 + \mu'_1 -1$ and $j < \beta_1 +
\mu_1 -1$, then $p$ does not divide $h_{\lambda}(i,j)$.  We assume that neither $\beta$ nor $\tau$ is empty, as simplified versions of the arguments that follow apply in those cases.  If $\tau'_1
+1 \leq i \leq \tau'_1 + \mu'_1 -1$ and $\beta_1 + 1 \leq j \leq
\beta_1 + \mu_1 -1$, then $h_{\lambda}(i,j)$ is a hook length of the
original $\mu$ and thus is not divisible by $p$ by assumption.  If $j \leq \beta_1$ and $i \leq \tau'_1$, so that the node $(i,j)$ is in the region of $\lambda$ above $\beta$ and to the left of $\tau$, note that 
\begin{equation} \label{decomp1}
h_{\lambda}(i,j) = h_{\lambda}(\tau'_1 + \mu'_1 ,j) + h_{\lambda}(i, \beta_1 + \mu_1) + h_{\lambda}(\tau'_1, \beta_1) - \tau_{\tau'_1} -
\beta'_{\beta_1}.
\end{equation} 
Here, $\tau_{\tau'_1}$ and $\beta'_{\beta_1}$ are the lengths of the last row of $\tau$ and the last column of $\beta$, respectively, and thus divisible by
$p$.  Since the first two terms and the last two terms on the right side of (\ref{decomp1}) sum are divisible by $p$, and the third is not, it follows that
$h_{\lambda}(i,j)$ is not divisible by $p$.  

If $j \leq \beta_1$ and
$\tau'_1 + 1 \leq i \leq \tau'_1 + \mu'_1 - 1$, so that $(i,j)$ is in the region of $\lambda$ directly to the left of $\mu$, including
the first column of $\mu$, then 
\begin{equation} \label{decomp2}
h_{\lambda}(i,j) = h_{\lambda}(\tau'_1 + \mu'_1 - 1, j) + h_{\mu}(x,y), 
\end{equation}
where $(x,y)$ is the node of $\mu$ corresponding to the node $(i, \beta_1)$
in $\lambda$.  Since $p$ divides the the first term on the right of (\ref{decomp2}), but not the second term, it follows that $h_{\lambda}(i,j)$ is not divisible by $p$.  A similar argument applies to the region of $\lambda$ directly
above $\mu$, including the first row of $\mu$.  Thus the only nodes of
$\lambda$ with a hook length divisible by $p$ are contained in
$\tau$ or $\beta$, and therefore $\lambda$ must be $p$-irreducible,
since $\beta$ and $\tau$ are $p$-irreducible.

Now we assume that $\lambda$ is $p$-irreducible and give the decomposition of $\lambda$.  If $\lambda$ is $p$-hook free, then we set $\lambda = \mu$, while $\tau$ and $\beta$ are empty, and we are done.  If $\lambda$ has a hook length divisible by $p$ at some node, then since $\lambda$ is $p$-irreducible, the entire row or entire column of that node must also consist of $p$-hooks.  So, we let $a$ be the smallest integer so that $h_{\lambda}(a,j)$ is divisible by $p$ for every node $(a,j)$ in row $a$, and we let $b$ be the smallest integer so that $h_{\lambda}(i,b)$ is divisible by $p$ for every node $(i,b)$ in column $b$.  It is possible that there is either no row of $p$-hooks or no column of $p$-hooks, in which case either $\beta$ or $\tau$ is empty, respectively, and simplified versions of the arguments which follow will apply.  So we assume that $a$ and $b$ are positive integers.

Let $(t,1)$ be any node in the first column of $\lambda$
with hook length divisible by $p$, and let $(1, u)$ be any node in the
first row of $\lambda$ with hook length divisible by $p$.  Note that $a \leq t$ and $b \leq u$.  We claim that $(t,u)$ is not a node of $\lambda$.  If $(t,u)$ is a node of $\lambda$, then 
$$h_{\lambda}(1, 1) = h_{\lambda}(t,1) + h_{\lambda}(1,u) - h_{\lambda}(t,u).$$
Since $(t,u)$ is in the same row as $(t,1)$, the $p$-irreducibility of $\lambda$ forces $h_{\lambda}(t,u)$ to be divisible by $p$, meaning that $h_{\lambda}(1,1)$ is divisible by $p$, a contradiction.  It follows that the row consisting of the nodes $(a,j)$ and the column consisting of the nodes $(i,b)$ cannot intersect. 

We now define the partitions $\tau$ and $\beta$.  By choice of $a$, for each $j$, all of the nodes in column $j$ of $\lambda$ above $(a,j)$ have hook lengths not divisible by $p$, and all of the nodes in row $a$ of $\lambda$ have
hook lengths divisible by exactly the same power of $p$ as
$h_{\lambda}(a,1)$.  Let $\beta$ be the partition defined by $\beta_r
= \lambda_{a+r-1}$ for $1 \le r \le \ell (\lambda)$.  Then $\beta$ is
necessarily a $p$-irreducible bottom.  Similarly, for each $i$, all of the nodes in row $i$ of $\lambda$ to the left of $(i,b)$ have hook lengths not divisible by $p$, and all of the nodes in column $b$ of $\lambda$ have hook lengths divisible by exactly the same power of $p$ as $h_{\lambda}(1,b)$.  Let $\tau$ be defined as the conjugate of the partition with parts $\tau'_s = \lambda'_{b+s-1}$ for $1 \le s \le \ell(\lambda')$.  Then $\tau$ must be a $p$-irreducible top.

Finally, we define the partition $\mu$ and show that it is $p$-hook free.  Recall that
$(a,1)$ is the first node in column 1 with hook length divisible by
$p$, and $(1, b)$ is the first node in row 1 with hook length
divisible by $p$.  Let $\mu$ be the partition consisting of all of the nodes of $\lambda$ in the set 
$$\{ (x,y) \, \big| \, \lambda'_b \le x < a \mbox{ and } \lambda_a \le y < b \}.$$
Thus $\mu$ is precisely the partition obtained from $\lambda$ by
removing all of the columns before and including $\beta_1 - 1$ and after and including column $b$, and all of the
rows before and including $\tau'_1 - 1$ and after and including row $a$.  Since row $a$ and column $b$ of $\lambda$ do not intersect in
$\lambda$, this process results in a well-defined partition $\mu$.  Finally, if any node $(x, y)$ in
$\mu$ had hook length divisible by $p$, then by the $p$-irreducibilty
of $\lambda$, either $(x, 1)$ or $(1, y)$ in $\lambda$ has hook length divisible by $p$.  But $(x,y)$ is above row $a$ or to the left of column $b$, which contradicts our choice of $a$ and $b$, and we are done.  By the choice of $a$ and $b$, it is clear that this decomposition of $\lambda$ is unique.
\end{proof}

One useful feature of the decomposition given in Theorem \ref{construction} is that $p$-regular and $p$-restricted partitions have a decomposition which is easy to describe.

\begin{cor} \label{decomppregpres}
Let $\lambda$ be a $p$-irreducible partition with decomposition $\lambda = \oplus(\tau, \mu, \beta)$, as given in Theorem \ref{construction}.  Then $\lambda$ is $p$-regular if and only if $\beta$ is empty, and $\lambda$ is $p$-restricted if and only if $\tau$ is empty.
\end{cor}
\begin{proof}
The second conclusion is equivalent to the first, since the conjugate of a $p$-regular partition is $p$-restricted.  Assume that $\lambda$ is not $p$-regular.  Then $\lambda$ has some part which occurs at least $p$ times.  If we consider the hook lengths of the nodes at the ends of these rows, then we must have a hook length $p$, since the hook lengths will increase from $1$ at the bottom, and there are at least $p$ nodes at the edge going up.  Since $\lambda$ is $p$-irreducible, and this column does not consist of all hook lengths divisible by $p$, then we must have a row of all hook lengths divisible by $p$.  This is equivalent to $\beta$ being nonempty.

Conversely, if $\beta$ is nonempty, then $\lambda$ has a row consisting of hook lengths all of which are divisible by $p$.  If we consider the last node in this row, we must have at least $p$ nodes including it and below.  But this means that this row is repeated at least $p$ times, and so $\lambda$ is not $p$-regular.\end{proof}

\section{The $p$-core of a $p$-irreducible partition}

In this section we will determine the structure of the $p$-core of a $p$-irreducible top or bottom.  We will then use this characterization, along with the decomposition of a $p$-irreducible partition given in Theorem \ref{construction}, to describe the $p$-core of a $p$-irreducible partition.  We subsequently develop a bijection between $p$-irreducible tops and certain ``smaller'' partitions, which will allow us to complete the counting argument in the next section.  For a partition $\lambda$, and a fixed prime $p$, we denote the $p$-core of $\lambda$ by $\hat{\lambda}$.

A {\it $p$-strip} is a horizontal or vertical strip of exactly $p$ boxes.  In the next lemma, we note that the $p$-core of any $p$-top, or $p$-bottom, respectively, can be obtained by removing $p$-strips from each row, or column, respectively.

\begin{lem} \label{easytauremoval}
If $\tau$ is a $p$-top with $k$ rows, then the $p$-core of $\tau$ is the partition $\hat{\tau}$ with parts $\hat{\tau}_i = (k-i)(p-1)$, $1 \leq i \leq k -1$, obtained by removing horizontal $p$-strips from $\tau$.  Similarly if $\beta$ is a $p$-bottom with $l$ columns, then the $p$-core of $\beta$, is the partition $\hat{\beta}$, the conjugate of the partition with parts $\hat{\beta}_i = (l -i)(p - 1)$, $1 \leq i \leq l - 1$, obtained by removing vertical $p$-strips from $\beta$.
\end{lem}

\begin{proof}
Note that the first and second statements are equivalent since a $p$-top is the conjugate of a $p$-bottom.  So assume that $\tau$ is a $p$-top with $k$ rows, and we prove the first statement by induction on $k$.  For $k = 1$ the statement is clear, since the $p$-core of $\tau$ is the empty partition. 

Now suppose $h_{\tau}(1,1) = rp$.  Then $\tau_1 = rp - k + 1$.  By the induction hypothesis, after taking the $p$-core of the bottom $k-1$ rows, we are left with the first row, and then the $k-2$ parts $(k-2)(p-1), (k-3)(p-1), \ldots, 2(p-1), p-1$.  Let us call this partition $\tilde{\tau}$.  Then 
$$\tilde{\tau}_1 - \tilde{\tau}_2 = rp - k + 1 - (k-2)(p-1) = p(r-k+1) + p-1.$$
Thus we may remove $r-k+1$ horizontal $p$-strips from the first row of $\tilde{\tau}$, leaving a row of length $(k-1)(p-1)$.  The lemma follows by noticing that the hook lengths of the first row of the resulting partition start with $(k-1)p - 1$, and decrease to $1$, while skipping multiples of $p$.   
\end{proof} 

\noindent
Lemma \ref{easytauremoval}, along with Theorem \ref{construction}, allows us to describe the $p$-core of any $p$-irreducible partition.  We first introduce a gluing of partitions which is a slightly modified version of the gluing in Definition \ref{decomp}.

\begin{defn} \label{coreglue}
Let $\tau$, $\mu$, and $\beta$ be three partitions, where $\tau$ and $\beta$ may be empty, but $\mu$ may only be empty if at least one of $\tau$ or $\beta$ is empty.  If $\tau = (\tau_1, \ldots, \tau_k)$, $\mu = (\mu_1, \ldots, \mu_l)$, and $\beta = (\beta_1, \ldots, \beta_m)$, then define the partition $\widehat{\oplus}(\tau, \mu, \beta)$ to be the partition with parts 
$(\tau_1 + \mu_1 + \beta_1, \, \tau_2 + \mu_1 + \beta_1, \ldots, \tau_k + \mu_1 + \beta_1, \, \mu_1 + \beta_1, \, \mu_2 + \beta_1, \ldots, \mu_l + \beta_1, \, \beta_1, \ldots, \beta_m)$.
\end{defn}

The difference between the pictures of $\oplus(\tau, \mu, \beta)$ and $\widehat{\oplus}(\tau, \mu, \beta)$ is that in the latter, the partitions $\tau$ and $\mu$, and $\beta$ and $\mu$, are glued at their corners, rather than in the manner of gluing described in Section \ref{decomposition}.

Recall the definition of the $p$-residual of a partition, given in Definition \ref{residual}.

\begin{cor} \label{lambdacore}
Assume that $\lambda$ is a $p$-irreducible partition which can be decomposed as $\oplus(\tau, \mu, \beta)$, as in Theorem \ref{construction}.  Then the $p$-core of $\lambda$ is the result of gluing, in the sense of Definition \ref{coreglue}, the $p$-cores of $\tau$ and $\beta$ to $\mu$.  That is, $\hat{\lambda} = \widehat{\oplus}(\hat{\tau}, \mu, \hat{\beta})$.

If $\tau$ has $k$ rows and $\beta$ has $l$ columns, then the $p$-residual $(t,b)$ of $\hat{\lambda}$ satisfies $k \leq t$ and $l \leq b$.
\end{cor}

\begin{proof}
Applying Lemma \ref{easytauremoval}, we may remove horizontal and vertical $p$-strips from $\lambda$, and we are left with $\nu = \widehat{\oplus}(\hat{\tau}, \mu, \hat{\beta})$.  We need only show that this partition is now $p$-hook free.  From the proof of Theorem \ref{construction}, we know that all of the hook-lengths in $\lambda$ for nodes above $\beta$ and to the left of $\tau$ are not divisible by $p$.  That is, for
$$ i \leq \ell(\tau) + \ell(\mu) - 1 \;\; \text{ and } \; \; j \leq \ell(\beta') + \ell(\mu') -1, $$
we have $v_p(h_{\lambda}(i,j)) = 0$.  In constructing $\nu$, we have only removed horizontal and vertical $p$-strips, and so these hook-lengths have changed by a multiple of $p$.  It follows that $\nu$ is $p$-hook free, and $\nu = \hat{\lambda}$.

The second statement follows immediately from the structure of $\hat{\tau}$ and $\hat{\beta}$ given in Lemma \ref{easytauremoval} and the definition of $p$-residual.
\end{proof}

The following lemma is an immediate consequence of the definition of hook-length.

\begin{lem}\label{lem:hookLength1}
Let $\sigma$ be a partition, and let $(i,j)$ be a node of $\sigma$ that is at the bottom of a column.  Then $h_\sigma (i,j) = \sigma_i + 1 - j$.
\end{lem}

The next lemma suggests a particularly useful correspondence between $p$-irreducible tops and smaller partitions.

\begin{lem} \label{shrinking}
Let $\nu = (\nu_1, \nu_2, \ldots, \nu_{k-1})$ be the partition given by $\nu_i = (k-i)(p-1)$ for $1 \le i \le k -1$, and let $\sigma = (\sigma_1, \ldots, \sigma_{k})$ be any partition with $\ell (\sigma) \le k$, so that $\sigma_i = 0$ for $\ell(\sigma) < i \leq k$.  Define a partition $\tau = (\tau_1, \ldots, \tau_{k})$ by $\tau_i = \nu_i + p \sigma_i$ for $i < k$ and $\tau_{k} = p \sigma_{k}$.  Then the following properties hold.
\begin{enumerate}
\item For every hook-length $q$ occuring in row $i$ of $\sigma$, there is a corresponding node in row $i$ of $\tau$ with a hook-length of $pq$, and the hook-lengths of $\tau$ obtained in this manner are the only $p$-hooks of $\tau$.
\item If the node in position $(i,j)$ of $\sigma$ is at the bottom of a column and has hook-length $q$, then the corresponding node of hook-length $pq$ in $\tau$ is in position $(i,\nu_i+p \sigma_i + 1 - pq)$.
\item If the hook-lengths $q$ and $q'$ occur in the same column of $\sigma$, then the corresponding hook-lengths $pq$ and $pq'$ of $\tau$ occur in the same column of $\tau$.
\end{enumerate}
\end{lem}

\begin{proof}
The proof is by induction on $|\sigma|$.  If $|\sigma| = 0$, then the result is obviously true since there are no $p$-hooks in either $\sigma$ or $\nu$ in this case.

Assume that $|\sigma| >0$, and let $m = \ell(\sigma)$ and $\tilde \sigma = (\tilde \sigma_1, \ldots, \tilde \sigma_{k})$, where
$$\tilde \sigma_i = \left\{ \begin{array}{ll} \sigma_i & \mbox{if $i \neq m$} \\ \sigma_i - 1 & \mbox{if $i = m$.} \end{array} \right.$$
The proof is straightforward if $m = 1$, so it is no loss to suppose that $m>1$.  Let $\tilde \tau = ( \tilde \tau_1, \ldots, \tilde \tau_{k})$ be the partition constructed from $\nu$ and $\tilde \sigma$ as in the statement of the lemma.  By induction, the statement of the lemma describes the relationship between $\tilde \tau$ and $\tilde \sigma$.

If $(i,j)$ is a node of $\tilde \sigma$, note that
$$h_\sigma (i,j) = \left\{ \begin{array}{ll} h_{\tilde \sigma} (i,j) & \mbox{if $i<m$ and $j \neq \sigma_m$} \\ h_{\tilde \sigma} (i,j) + 1 & \mbox{if $i<m$ and  $j = \sigma_m$} \\ h_{\tilde \sigma} (i,j) + 1 & \mbox{if $i=m$.} \end{array} \right.$$
A similar relationship exists between the hook-lengths of $\tau$ and $\tilde \tau$.  Let $I = \{ j \in \z \, | \, \nu_m + p \tilde \sigma_m + 1 \le j \le \nu_m + p \sigma_m \}$, and observe that
\begin{equation} \label{tauhooks} h_\tau (r,s) = \left\{ \begin{array}{ll} h_{\tilde \tau} (r,s) & \mbox{if $r \neq m$ and  $s \not\in I$} \\ h_{\tilde \tau} (r,s) + 1 & \mbox{if $r \neq m$ and  $s \in I$} \\ h_{\tilde \tau} (r,s) + p & \mbox{if $r = m$.} \end{array} \right.
\end{equation}
whenever $(r,s)$ is a node of $\tilde \tau$.  Thus the hook-lengths of $\sigma$ that differ from those of $\tilde \sigma$ are in row $m$ or in column $\sigma_m$, and the hook-lengths of $\tau$ that differ from those of $\tilde \tau$ are in row $m$ or in a column with index in $I$.

If $j< \sigma_m$ and $(m,j)$ is a node of $\sigma$ with hook-length $q$, then the node $(m,j)$ of $\tilde \sigma$ has a hook-length of $q-1$.  By induction, there is a corresponding hook-length of $p(q-1)$ in row $m$ and column $\nu_m + p \tilde \sigma_m + 1 - p(q-1)$ of $\tilde \tau$.  Since the hook-length for the node $(m, \nu_m + p \tilde \sigma_m + 1 - p(q-1))$ of $\tau$ is $p(q-1) + p = pq$, and since $\nu_m + p \tilde \sigma_m + 1 - p(q-1) = \nu_m + p \sigma_m + 1 - pq$, we see that $\tau$ has a hook-length of $pq$ that corresponds to the hook-length of $q$ in $\sigma$, and it is indeed in the claimed position.  Note that outside of $I$, the hook lengths of $\tau$ are the same as the hook lengths of $\tilde \tau$, except that in row $m$ the hook lengths of $\tau$ are $p$ greater than the hook lengths of $\tilde \tau$.

It remains to show that column $\sigma_m$ of $\sigma$ corresponds to column $\nu_m+p \sigma_m +1 - ph_\sigma (m, \sigma_m)$ of $\tau$ and that none of the nodes of $\tau$ in other columns in $I$ contain $p$-hooks.

By induction, the $p$-hooks of $\tilde \tau$ occur in columns, and the node in position $(m-1, \sigma_m)$ of $\tilde \sigma$ corresponds to the node of $\tilde \tau$ in position
$$(m-1, \nu_{m-1} + p \tilde \sigma_{m-1}+1-ph_{\tilde \sigma}(m-1,\sigma_m)),$$
From Lemma \ref{lem:hookLength1}, we see that
$$\nu_{m-1} + p \tilde \sigma_{m-1}+1-ph_{\tilde \sigma}(m-1,\sigma_m) = \nu_m +p \sigma_m \in I.$$
This implies that column $\sigma_m$ of $\tilde \sigma$ corresponds to the column $\nu_m +p \sigma_m \in I$ of $\tilde \tau$, and thus by induction $h_{\tilde \tau}(r, \nu_m+ p \sigma_m) \equiv 0 \quad (\mbox{mod } p)$ for $1 \le r \le m-1$.  If $s-1 \in I$, then $h_{\tilde \tau}(r,s-1) = 1 + h_{\tilde \tau}(r,s)$, and we may generalize that 
\begin{equation}\label{eqn:hooksInI}
h_{\tilde \tau}(r, \nu_m + p \sigma_m - t) \equiv t \quad (\mbox{mod } p)
\end{equation}
for $0 \le t \le p-1$.  

As $h_\sigma (m,\sigma_m) = 1$, it follows that $\nu_m+p \sigma_m +1 - ph_\sigma (m, \sigma_m) = \nu_m + p \tilde \sigma_m + 1$.  Note also that $h_\tau (m, \nu_m+ p \tilde \sigma_m + 1) = p$ and $h_\tau (m,s) \not\equiv 0 \quad (\mbox{mod } p)$ if $\nu_m+ p \tilde \sigma_m + 1 \neq s \in I$.  For $1 \le i \le m-1$, let $q_i = h_\sigma (i, \sigma_m)$ and $\tilde q_i = h_{\tilde \sigma} (i, \sigma_m)$, and let $q_m = h_\sigma (m, \sigma_m) = 1$.  By induction, we know that $h_{\tilde \tau}(i, \nu_m+p \sigma_m) = p \tilde q_i$.  Since node $(m, \sigma_m)$ of $\sigma$ is at the bottom of column $\sigma_m$, conditions (2) and (3) are equivalent to the condition that
$$h_\tau (i, \nu_m + p \tilde \sigma_m + 1) = pq_i$$
for $1 \le i \le m$.  We have seen that $h_\tau ( m, \nu_m + p \tilde \sigma_m + 1) = p = p q_m$.  For $i<m$, note that $q_i = \tilde q_i+1$, and thus we have 
\begin{align*}
h_\tau (i, \nu_m + p \tilde \sigma_m + 1) &= 1 + h_{\tilde \tau}(i, \nu_m + p \tilde \sigma_m + 1) \\
&= 1 + h_{\tilde \tau}(i, \nu_m + p \sigma_m) + p-1 \\
&= 1+ p \tilde q_i + p-1 \\
&= pq_i.
\end{align*}
Because $h_\tau (r,s) = h_{\tilde \tau} (r,s) + 1$ for $s \in I$, it follows from (\ref{eqn:hooksInI}) that
$$h_\tau (r, \nu_m + p \sigma_m - t) = h_{\tilde \tau} (r, \nu_m + p \sigma_m - t) + 1 \equiv t+1 \quad (\mbox{mod } p)$$
for $1 \le r \le m-1$ and $0 \le t \le p-1$.  This implies that the only $p$-hooks occuring in some column indexed by $I$ must occur in column $\nu_m + p \tilde \sigma_m + 1$.  
\end{proof}

\begin{rem}\label{rem:pRegular}
Lemma \ref{shrinking} may be interpreted as stating that the
$p$-hooks of the partition $\tau$ occur in columns.
\end{rem}

The following consequence of Lemma \ref{shrinking} is a key ingredient to the proof of the main counting theorem in the next section.

\begin{cor} \label{correspondence}
Let $\nu$, $\sigma$, and $\tau$ be as in Lemma \ref{shrinking}.  Then $\tau$ is $p$-irreducible if and only if $\sigma$ is $p$-irreducible and $p$-regular.  In the case that $\sigma$ is $p$-irreducible and $p$-regular, $\tau$ is a $p$-irreducible top if and only if $\ell (\sigma ) = k$. 
\end{cor}
\begin{proof}
Assume that $\tau$ is $p$-irreducible.  Suppose that $v_p(h_{\sigma}(i,j))> 0$ for some node $(i,j)$ of $\sigma$, and let $(i',j)$ be some other node of $\sigma$ in the same column.  By Lemma \ref{shrinking}, there are corresponding nodes $(i, l)$ and $(i', l)$ of $\tau$ such that 
\begin{equation}\label{eqn:correspnd1}
h_\tau (i, l) = p h_\sigma (i, j) \quad \mbox{and} \quad h_\tau
(i', l) = p h_\sigma (i',j).
\end{equation}
Since the $p$-hooks of $\tau$ occur in columns and $\tau$ is by assumption $p$-irreducible, it follows from (\ref{eqn:correspnd1}) that $v_p(h_{\sigma}(i,j)) = v_p(h_{\sigma}(i',j))$, and thus $\sigma$ is $p$-irreducible.  If $\sigma$ were not $p$-regular, there would necessarily be some column of $\sigma$ containing both a $p$-hook and non-$p$-hooks.  As this is not true of $\sigma$, it must be that $\sigma$ is $p$-regular.

Suppose that $\sigma$ is $p$-irreducible and $p$-regular.  By Corollary \ref{decomppregpres}, the $p$-hooks of $\sigma$ occur in columns, and thus 
\begin{equation}\label{eqn:correspond2}
v_p(h_{\sigma}(i,j)) = v_p(h_{\sigma}(i',j)) \quad \mbox{for all $i,i',j$}.
\end{equation}
The $p$-hooks of $\tau$ must also occur in columns as described in Lemma \ref{shrinking}.  If column $l$ of $\tau$ contains $p$-hooks, then Lemma \ref{shrinking} implies that there exists $j$ such that 
$$h_\tau (i, l) = p h_\sigma (i, j) \quad \mbox{and} \quad h_\tau (i', l ) = p h_\sigma (i',j)$$
for all (rows) $i,i'$.  Hence, (\ref{eqn:correspond2}) yields that 
$$v_p(h_\tau (i,l)) = v_p (h_\tau (i',l)),$$
and it follows that $\tau$ is $p$-irreducible. 

Finally, if $\tau$ is $p$-irreducible, then $\tau$ is a
$p$-irreducible top if and only if the left-most column of $\tau$
consists entirely of $p$-hooks.  From Lemma \ref{shrinking}, this is
clearly the case if and only if $\sigma_{k} \neq 0$.
\end{proof} 

\begin{rem} \label{rem:conj}
It is clear that Lemma \ref{shrinking} and Corollary \ref{correspondence} can be modified by conjugation, so that we may also apply the results to $p$-bottom partitions.
\end{rem}

\section{Main Results} \label{mainthm}

\subsection{The case $p>2$}  We begin with a technical lemma about the $p$-residual of a partition when $p>2$.

\begin{lem} \label{reslemma}
Let $p>2$ and $\nu$ be any partition with $p$-residual $(t,b)$.  Then we have
$$ t + b \leq \frac{\ell(\nu) + \ell(\nu')}{p} + 2. $$
Equality occurs if and only if either $(t,b) = (1,\ell(\nu')+1)$, or $(t,b) = (\ell(\nu) + 1, 1)$, or $t>1$ and $b>1$ and $(t-1,b-1)$ is a node in $\nu$, but $(t,b)$ is not.
\end{lem}
\begin{proof}
Suppose that $t>1$ and $b>1$.  From Definition \ref{residual}, we have $\nu_{i} - \nu_{i + 1} = p - 1$ for $0<i<t$.  Then for each $i<t$, if $k$ is such that $\nu_{i+1} - (p -1) < k < \nu_{i+1}$, we have $\nu'_k - \nu'_{k + 1} = 0$, and for $k  = \nu_{i+1}$, we have $\nu'_k - \nu'_{k+1} = 1 < p -1$, since $p > 2$.  So if for some $j$ we have $\nu'_j - \nu'_{j + 1} = p - 1$, then $j$ must satisfy $j \leq \nu_{i+1} - (p-1)$ for every $i<t$.  Since the $p$-residual of $\nu$ is $(t,b)$, then we must have $b - 1 \leq \nu_{t-1} - (p-1)$.  Using the fact that $\nu_i = \nu_{i - 1} - (p -1)$ for $i<t$, we obtain,
$$ \nu_{t-1} - (p-1) = \nu_1 - (t-1)(p-1) = \ell(\nu') - (t-1)(p-1), $$
since $\nu_1 = \ell(\nu')$.  So now we have
\begin{equation} \label{rows} 
b - 1 \leq \ell(\nu') - (t-1)(p-1).
\end{equation} 
By the exact same argument, but interchanging the roles of rows and columns, we have
\begin{equation} \label{cols}
t - 1 \leq \ell(\nu) - (b-1)(p-1).
\end{equation}
It is easy to see that (\ref{rows}) and (\ref{cols}) also hold when $t=1$ or $b=1$.  Adding (\ref{rows}) and (\ref{cols}) and simplifying yields
$$ t + b \leq \frac{\ell(\nu) + \ell(\nu')}{p} + 2. $$
Equality occurs when we have equality in (\ref{rows}) and (\ref{cols}).  If $t=1$, the maximum value for $b$ is $\ell(\nu') + 1$, and similarly if $b=1$, the maximum value for $t$ is $\ell(\nu) + 1$, and equality occurs in both of these cases.  Finally, if $b>1$ and $t>1$, (\ref{rows}) and (\ref{cols}) imply that column $b - 1$ and row $t - 1$ of $\nu$ intersect, while column $b$ and row $t$ do not, which is equivalent to $(t - 1, b - 1)$ being a node in $\nu$ while $(t,b)$ is not.
\end{proof}

We are now able to describe how to count and construct the $p$-irreducible Specht modules in a given block, for $p > 2$.

\begin{thm} \label{COVproof}
Let $B$ be a $p$-block of $S_n$ (with $p > 2$), corresponding to the $p$-hook free partition $\nu$, of weight $w$ and $p$-residual $(t,b)$.  Then the number of $p$-irreducible Specht modules in $B$ is equal to the number of ordered pairs of $p$-regular $p$-irreducible partitions $(\alpha, \gamma)$, such that $|\alpha| + |\gamma| = w$, $\ell(\alpha) \leq t$, $\ell(\gamma) \leq b$, and $\ell(\alpha) + \ell(\gamma) \leq t+b-1$ if $t + b = \frac{\ell(\nu)+ \ell(\nu')}{p} + 2$.
\end{thm}
\begin{proof} First recall Nakayama's conjecture, which states that the Specht modules $S^{\lambda}$ and $S^{\rho}$ of $S_n$ are in the same $p$-block if and only if $\lambda$ and $\rho$ have the same $p$-core.  So, we describe all $p$-irreducible $\lambda$ such that $\hat{\lambda} = \nu$.  

Suppose that $\lambda$ is a $p$-irreducible partition of $n$ such that $\hat{\lambda} = \nu$.  From Lemma \ref{easytauremoval} and Corollary \ref{lambdacore}, we know that $\nu$ is obtained from $\lambda$ by removing a total of exactly $w$ horizontal and vertical $p$-strips.  So, to obtain any $p$-irreducible $\lambda$ such that $\hat{\lambda} = \nu$, we must add a total of $w$ horizontal and vertical $p$-strips to $\nu$.  If we add horizontal $p$-strips to $k$ rows and vertical $p$-strips to $l$ columns, then it follows from Corollary \ref{lambdacore} that we must have $k \leq t$ and $l \leq b$, where $(t,b)$ is the $p$-residual of $\nu$.

It follows from the structure of $\nu$ given in Lemma \ref{easytauremoval} and Corollary \ref{lambdacore} that if we are adding horizontal $p$-strips to the first $t$ rows of $\nu$, then in order for the result to be a partition, then we can add at most as many $p$-strips to each row as we did the row immediately above.  This is also true for adding vertical $p$-strips to the first $b$ columns of $\nu$.  Thus we may view the set of horizontal $p$-strips added and vertical $p$-strips added as partitions, $\alpha$ and $\gamma$, where the parts of $\alpha$ and $\gamma$ are given by the number of $p$-strips added to each successive row and column of $\nu$, respectively.  Then we must have $|\alpha| + |\gamma| = w$, and from the comments in the above paragraph, $\ell(\alpha) = k \leq t$ and $\ell(\gamma) = l \leq b$.

The partitions $\alpha$ and $\gamma$ must be such that the resulting partition $\lambda$ obtained is $p$-irreducible.  It follows from Corollary \ref{correspondence} that this is satisfied exactly when $\alpha$ and $\gamma$ are both $p$-regular and $p$-irreducible.

Finally, consider the case when $t + b$ is maximal, which occurs in the situations described in Lemma \ref{reslemma}.  In each of these cases, it means that when putting $p$-strips on $\nu$ to obtain $p$-irreducibles in $B$, there is exactly one place where we could put a horizontal $p$-strip or a vertical $p$-strip, which is in position $(t,b)$ in the case that $t>1$ and $b>1$.  Since we may only put one or the other in this position, we must restrict the sum of the lengths of $\alpha$ and $\gamma$ to be 1 less than the maximum, so that we must have $\ell(\alpha) + \ell(\gamma) \leq t+b-1$ in the case that $t+b = \frac{\ell(\nu)+ \ell(\nu')}{p} + 2$.
\end{proof}

Given an arbitrary block of $S_n$, one might expect that the $p$-residual of the block is a pair of small numbers.  For the case of the smallest possible $p$-residual, we can say precisely how many $p$-irreducible Specht modules are in the block.

\begin{cor} \label{smallres}
Let $p>2$, and let $B$ be a block of $S_n$, corresponding to the $p$-hook free partition $\nu$, of weight $w$.  Suppose that the $p$-residual of $B$ is $(1,1)$, so that either $\nu$ is empty or $\nu$ is nonempty and satisfies $\nu_1 - \nu_2 \neq p-1$ and $\nu'_1 - \nu'_2 \neq p-1$.  Then the the number of $p$-irreducible Specht modules in $B$ is equal to $2$ if $\nu$ is empty, and is equal to $w+1$ otherwise.
\end{cor}
\begin{proof} If $\nu$ is empty, then from Theorem \ref{COVproof}, we may either add $w$ horizonatal $p$-strips, or $w$ vertical $p$-strips, but not both.  So there are only $2$ possible $p$-irreducible Specht modules in this case.  Otherwise, from Theorem \ref{COVproof}, the number of $p$-irreducible Specht modules in $B$ is the number of ordered pairs of $p$-irreducible and $p$-regular partitions with at most one part, the sum of whose sizes is $w$.  Since every partition with at most one part is $p$-irreducible and $p$-regular, then we can just count the number of ordered pairs of non-negative numbers with sum $w$, which is $w+1$.   
\end{proof}

\subsection{The case $p=2$} \label{twocase}
The description of $2$-irreducible Specht modules of $S_n$ is completed in \cite{JM}, and the situation is a bit different than in the case $p \neq 2$.  We continue saying that the Specht module $S^{\lambda}$ defined over a characteristic $0$ field (and thus over the integers) is $2$-irreducible if it remains irreducible when reduced modulo $2$, which is equivalent to the Specht module $S^{\lambda}_{\mathbb{F}_2}$ being irreducible.

We say that a partition $\lambda$ is $2$-irreducible if it satisfies the condition in Theorem \ref{fayers}.  That is, we say that $\lambda$ is $2$-irreducible when $\lambda$ does not have nodes $(i,j)$, $(i,y)$, and $(x,j)$ such that 
$$v_2(h_{\lambda}(i,j)) > 0, \; v_2(h_{\lambda}(i,y)) \neq v_2(h_{\lambda}(i,j)), \text{ and } v_2(h_{\lambda}(x,j)) \neq v_2(h_{\lambda}(i,j)).$$

Unlike in the case $p \neq 2$, it is not always the case that if $\lambda$ is $2$-irreducible then $S^{\lambda}$ is $2$-irreducible.  Rather, the situation is as follows.

\begin{thm}[James and Mathas, 1999] \label{twoirred}
If $n \neq 4$, the Specht module $S^{\lambda}$ of $S_n$ is $2$-irreducible if and only if $\lambda$ is either $2$-regular or $2$-restricted, and $\lambda$ is $2$-irreducible.  If $n=4$, then $S^{\lambda}$ is $2$-irreducible if and only if the previous conditions hold or if $\lambda = (2^2)$.
\end{thm} 

We have the following description of $2$-irreducible Specht modules which appear in a $2$-block of $S_n$.

\begin{thm} \label{twocount}
Let $n \neq 4$, and let $B$ be a $2$-block of $S_n$, corresponding to the $2$-hook free partition $\nu$, of weight $w$.  Then the number of $2$-irreducible Specht modules in $B$ is equal to twice the number of $2$-regular $2$-irreducible partitions $\alpha$ of $w$ such that $\ell(\alpha) \leq \ell(\nu) + 1$.
\end{thm}
\begin{proof} Let $\nu$ be the $2$-hook free partition corresponding to the block $B$.  The proof is very similar to the proof of Theorem \ref{COVproof}, except by Corollary \ref{decomppregpres} and Theorem \ref{twoirred}, we either add horizontal $2$-strips or vertical $2$-strips to $\nu$, but not both.  Also note that the $2$-residual $(t,b)$ of any nonempty partition $\nu$ satisfies $b= \ell(\nu) = \ell(\nu') = t$.  If $\alpha$ is the partition formed by the process of adding horizontal $2$-strips, and $\gamma$ is the corresponding partition from adding vertical $2$-strips, then either $\alpha$ or $\gamma$ is empty, while the other is a partition of $w$.  From Corollary \ref{correspondence}, if $\alpha$ or $\gamma$ is non-empty, then it must be $2$-irreducible and $2$-regular.  We may add $2$-strips to every row or column of $\nu$, and we may even add an extra row or column, from the structure of $2$-hook free partitions, so that the length of $\alpha$ or $\gamma$ may be one more than the number of rows of $\nu$, which is equal to the number of columns of $\nu$.  So $\alpha$ and $\gamma$ must have the exact same structure, and we can just count the certain type of partition twice.  
\end{proof}

In general, one may be interested in $p$-irreducible Specht modules which correspond to only $p$-regular or $p$-irreducible partitions.  In this situation, we have the following result, which is true for any prime $p$, keeping in mind that the $2$-residual of a $2$-block only depends on the length of the corresponding $2$-hook free partition.

\begin{cor} \label{pregpres}
Let $B$ be a $p$-block of $S_n$ (for any $p$), corresponding to the $p$-hook free partition $\nu$, of weight $w$ and $p$-residual $(t,b)$.  Then the number of $p$-regular $p$-irreducible Specht modules in $B$ is equal to the number of $p$-regular $p$-irreducible partitions $\alpha$ of $w$ such that $\ell(\alpha) \leq t$, and the number of $p$-restricted $p$-irreducible Specht modules in $B$ is equal to the number of $p$-regular $p$-irreducible partitions $\gamma$ of $w$ such that $\ell(\gamma) \leq b$.
\end{cor}

\begin{proof}  Note that by Corollary \ref{decomppregpres}, a
  $p$-irreducible partition $\lambda$ is $p$-regular if and only if
  $\lambda$ has no bottom part.  Thus the number of $p$-regular
  $p$-irreducible partitions in $B$ is equal to the number of $p$-regular
  $p$-irreducible partitions $\alpha$ with length at most $t$.  A
  similar argument holds for the $p$-restricted $p$-irreducible
  partitions in $B$.
\end{proof}

\section{An Example} \label{example}
In this section, we give an example of how the above methods can
be used to count and explicitly construct all of the irreducible
Specht modules in a given block.  Suppose that $p=5$ and the block $B$ is given by the $5$-hook free partition $\nu = (17, 13, 9, 5^2, 3^3, 2^4, 1^4)$, pictured in Figure \ref{ex}.

\begin{figure}
\caption{}
\yng(17,13,9,5,5,3,3,3,2,2,2,2,1,1,1,1)
\label{ex}
\end{figure}

Notice that $\nu$ is a partition of 70, and that $\nu_1-\nu_2 = \nu_2 - \nu_3 = \nu_3-\nu_4 = p-1 = 4$.  Therefore $t = 4$, and similarly $b= 3$.
Now suppose that $B$ is a block of weight 8, so that $\nu$ corresponds to a block of $S_{110}$.  Then by applying Theorem \ref{COVproof}, the partitions corresponding to $5$-irreducible Specht modules in $B$ are indexed by ordered
pairs of partitions $\alpha$ and $\gamma$ such that $\alpha$ has at
most 4 rows, $\gamma$ has at most $3$ rows, $|\alpha| + |\gamma| = 8$, and $\alpha$ and $\gamma$ are both $5$-irreducible and 5-regular.

We now only need to count the pairs of partitions $(\alpha, \gamma)$
with the above properties.  We let $p_{\alpha}(k)$ denote the number
of partitions of $k$ with at most four rows that are both
5-irreducible and 5-regular, and let $p_{\gamma}(k)$ denote the
number of partitions of $k$ with at most three rows that are both
5-irreducible and 5-regular, and we define $p_{\alpha}(0) = p_{\gamma}(0) = 1$.  It may be easily checked by hand that the sequence $(p_{\alpha}(0), p_{\alpha}(1), \ldots, p_{\alpha}(8)) = (1, 1, 2, 3, 5, 3, 6, 6, 8)$ and the sequence $(p_{\gamma}(0), p_{\gamma}(1), \ldots, p_{\gamma}(8))$ = $(1, 1, 2, 3, 4, 3, 5, 4, 5)$.  Thus the number of $p$-irreducible Specht modules in this block is given by
$$\sum_{k=0}^8 p_{\alpha}(k) p_{\gamma}(8-k) = 83.$$  Moreover, these
partitions can be easily constructed.  For instance, if $\alpha =              (2^2, 1)$ and $\gamma = (2,1)$, then the partition $\lambda$ corresponding to the 5-irreducible
Specht module $S^{\lambda}$ in $B$ is constructed by adding ten boxes to each of the first two rows of $\nu$, five boxes to the third row of $\nu$, ten boxes to the first column of $\nu$, and five boxes to the second column of $\nu$.  The resulting 5-irreducible partition is $\lambda = (27, 23, 14, 5^2, 3^2, 2^9, 1^9)$.

\end{document}